\documentclass[reqno]{amsart}

\usepackage[dvipsnames]{xcolor}

\usepackage{iftex}
\ifPDFTeX 
  \usepackage[utf8]{inputenc}
  \usepackage[T1]{fontenc}
\else 
  \usepackage{fontspec}
\fi
\usepackage{lmodern}

\usepackage{geometry}
\usepackage[textwidth=28mm]{todonotes}
\usepackage{mathtools}
\usepackage{yfonts}
\usepackage[foot]{amsaddr}

\usepackage{algorithm}
\usepackage{algpseudocode}

\usepackage[pdfencoding=auto]{hyperref}
\hypersetup{
    colorlinks,
    citecolor=green,
    filecolor=black,
    linkcolor=blue,
    urlcolor=black
}

\newtheorem{Definition}{Definition}

\newtheorem{Remark}{Remark}

\newcommand{\dt}{\,\partial_t\, }

\newcommand{\ddt}{\frac{\dd}{\dd t}}

\newcommand{\dd}{\,\mathrm{d}\,}

\newcommand{\N}{\mathbb{N}}

\newcommand{\R}{\mathbb{R}}

\newcommand{\M}{\mathcal{M}}

\newcommand{\iph}{{i+\frac{1}{2}}}
\newcommand{\imh}{{i-\frac{1}{2}}}
\newcommand{\jph}{{j+\frac{1}{2}}}
\newcommand{\jmh}{{j-\frac{1}{2}}}
\newcommand{\ep}{\varepsilon}

\title[Kinetic Fluid Parareal Method]{A Parareal in time numerical method for the collisional Vlasov equation in the hyperbolic scaling}
\author{Tino Laidin}
\address[Tino Laidin]{Univ Brest, CNRS UMR 6205, Laboratoire de Mathématiques de Bretagne Atlantique, F-29200 Brest, France}
\author{Thomas Rey}
\address[Thomas Rey]{Université Côte d’Azur, CNRS, LJAD, Parc Valrose, F-06108 Nice, France}
\email{tino.laidin@univ-brest.fr}
\email{thomas.rey@univ-cotedazur.fr}

\begin{document}

\begin{abstract}
    We present the design of a multiscale parareal method for kinetic equations in the fluid dynamic regime. The goal is to reduce the cost of a fully kinetic simulation using a parallel in time procedure. Using the multiscale property of kinetic models, the cheap, coarse propagator consists in a fluid solver and the fine (expensive) propagation is achieved through a kinetic solver for a collisional Vlasov equation. To validate our approach, we present simulations in the 1D in space, 3D in velocity settings over a wide range of initial data and kinetic regimes, showcasing the accuracy, efficiency, and the speedup capabilities of our method.
    \\[1em]

\end{abstract}

\maketitle

\section{Introduction}
The simulation of many complex phenomena involving interacting particles can be modelled through either fluid or kinetic descriptions. However, the validity of the models varies highly depending on the application. In particular, to accurately describe rarefied gases or charged particles in a device, a fluid description given by the Euler, Navier-Stokes or drift-diffusion equations may break down. Typically, it occurs around shocks or because of boundary layers and a kinetic description is therefore necessary. Nevertheless, the main drawback of the latter approach lies in its numerical cost. Indeed, deterministic kinetic simulations suffer even more strongly than usual solvers from the curse of dimensionality because of the large size of the phase space. As a consequence, one wants to resort to kinetic simulation sparingly.

In this work, we are interested in the approximation of solutions to the following scaled collisional Vlasov equation:
\begin{equation}\label{para:VcolScaled}
    \left\lbrace\begin{aligned}
        &\dt f^\ep(t,x,v) + \frac{1}{\ep^{\alpha}}\big(v\cdot\nabla_x\,f^\ep(t,x,v) - E\,\cdot\,\nabla_v\,f^\ep(t,x,v)\big)= \frac{1}{\ep^{\alpha+1}}\mathcal{Q}(f^\ep),\\
        &f(0,x,v) = f_0(x,v),
    \end{aligned}\right.
\end{equation}
where $f(t,x,v)$ is the distribution function of particles, $t\geq0$, $x\in\Omega_x\subset\R^{d_x}$, $v\in\R^{d_v}$. In position, the domain $\Omega_x$ is bounded, and boundary conditions will be presented later. This equation is used to model several phenomena including rarefied gases ($E=0$) or charged particles in a tokamak, or electronic device where $E$ typically is gradient of a potential that solves a Poisson equation. The parameter $\ep$ is the Knudsen number, which is the ratio between the mean free path of particles and the typical length scale of observation. Depending on the value of $\alpha$, the asymptotic model as $\ep$ tends to $0$ will either be given by a \emph{hydrodynamic} system ($\alpha=0$) \cite{BardosGolseLevermore1991,GolseSaintRaymond2001}, such as Euler's equations, or a \emph{diffusive} one ($\alpha=1$) \cite{Poupaud1991,PoupaudSoler2000}, typically a drift-diffusion equation.

A main component of \eqref{para:VcolScaled} is the Boltzmann-like collision operator $\mathcal{Q}$. In the following, we will assume that it satisfies the following classical properties:
\begin{itemize}
    \item It preserves mass, momentum and energy: \begin{equation}\label{ConsColl}
        \int_{\R^3} \mathcal{Q}(f)\begin{pmatrix} 1\\v\\\frac{|v|^2}{2}\end{pmatrix} \dd v = 0_{\R^5}.
    \end{equation}
    \item It dissipates the Boltzmann entropy: \begin{equation}\label{DissEntro}
        \int_{\R^3} \mathcal{Q}(f)\log(f) \dd v \leq 0.
    \end{equation}
    \item Its equilibria are given by Maxwellian distributions: \begin{equation}\label{EquiMaxwellienne}
        \mathcal{Q}(f)=0 \iff f=\M_{\rho,u,\theta}(v)\coloneqq\frac{\rho}{(2\pi\theta)^{d/2}}\exp\left(-\frac{|v-u|^2}{2\theta}\right), \quad \forall v \in \R^3,
    \end{equation}
\end{itemize}
where $\rho$, $u$ and $\theta$ denote the density, mean velocity and temperature respectively. In regard to property \eqref{EquiMaxwellienne} one obtains, at least formally \cite{jungel_transport_2009}, that in the limit $\ep\rightarrow0$, the distribution $f^\ep$ tends towards a Maxwellian distribution whose moments are solution to a fluid model. In the hydrodynamic scaling ($\alpha=0$), the moments solve Euler's equations for a monoatomic gas:
\begin{equation}\label{para:Euler}
    \left\lbrace\begin{aligned}
        & \dt \rho +\nabla_x\,\cdot(\rho u)=0,\\
        & \dt(\rho u)+\nabla_x\,\cdot(\rho u\otimes u) + \nabla_x(\mathbb{P})-\rho E=0,\\
        & \dt \mathcal{E} + \nabla_x\,\cdot((\mathcal{E}+\mathbb{P})u) - \rho u\cdot E=0,
    \end{aligned}\right.
\end{equation}
where $\mathbb{P}=\rho\theta I$ is the pressure and the energy $\mathcal{E}$ is related to the temperature through the relation,
\begin{equation*}
    \mathcal{E} = \frac{1}{2}\left(\rho|u|^2+d_v\rho\theta\right).
\end{equation*}
In the diffusive scaling ($\alpha=1$), they solve a drift-diffusion equation ($u=0$, $\theta=1$):
\begin{equation}\label{para:DDbis}
    \dt\rho-m_2\mathrm{div}_x J=0,\text{ where } J=\nabla_x\rho-E\rho,
\end{equation}
and $m_2$ is the a space averaged second moment.
A widely used collision operator that we consider in this work, and which is simpler than the quadratic Boltzmann operator, is the so-called BGK operator \cite{BhatnagarGrossKrook1954}. It satisfies the above properties and is a relaxation operator with Maxwellian equilibrium:
\begin{equation}\label{para:BGK}
    \mathcal{Q}(f)=\frac{\tau}{\varepsilon}(\M_{\rho,u,\theta}-f),
\end{equation}
where $\tau$ depends on the moments $\rho$ and $\theta$. This simple model gives the correct asymptotic models as $\ep\rightarrow0$, and an extended study using a Chapman-Enskog expansion \cite{Struchtrup2005book, jungel_transport_2009} shows that higher order systems can be obtained. For example, in the hydrodynamic scaling, one can obtain the Navier-Stokes equations:
\begin{equation}\label{para:CNS}
    \left\lbrace\begin{aligned}
        & \dt \rho +\nabla_x\,\cdot(\rho u)=0,\\
        & \dt(\rho u)+\nabla_x\,\cdot(\rho u\otimes u) + \nabla_x(\mathbb{P})-\rho E=\ep\nabla_x\,\cdot\,\left(\mu\sigma(u)\right),\\
        & \dt \mathcal{E} + \nabla_x\,\cdot((\mathcal{E}+\mathbb{P})u)- \rho u \cdot E=\ep\nabla_x\,\cdot\,\left(\mu\sigma(u)u+\kappa\nabla_x\,\theta\right),
    \end{aligned}\right.
\end{equation}
with the strain rate tensor $\sigma(u)$ given by
\begin{equation*}
    \sigma(u)=\left(\nabla_x\,u+(\nabla_x\,u)^T\right) - \frac{2}{3}(\nabla_x\,\cdot\,u)I.
\end{equation*}

Several approaches have been designed to significantly reduce the computational cost of kinetic solvers. An inspiration for this work are the so-called \emph{hybrid methods}. The idea behind these approaches is to achieve a coupling between a cheap, low dimensional, fluid model and the expensive, high dimensional, kinetic one. This type of techniques \cite{TiwariKlar1998,CrouseillesDegondLemou2004,FilbetRey2015} mainly relies on some kind of domain decomposition in position or velocity. To accurately describe the solution, domain indicators that assess the validity of the cheap fluid description over the expensive kinetic one are often needed \cite{LevermoreMorokoffNadiga1998,Tiwari1998,DegondDimarcoMieusens2007,KolobovArslanbekovAristovFrolovaZabelok2007}. In addition to the computation of subdomains (that can be costly depending on the physical regime considered), one must also deal with boundary conditions between different kinds of models.

The last type of methods aiming at reducing the cost of kinetic simulations that we would like to mention  are the so-called Dynamical Low Rank methods, that directly deal with the matrix of the discretized system \cite{KockLubich2007,EinkemmerHuWang2021} through some type of SVD-like decomposition. 

\medskip

\textbf{The parareal method.} Our goal in this work is to take a new strategy, which focuses on the time integration of the multiscale kinetic equation considered. More precisely, we design a \emph{Parallel in Time} (PinT) method that takes advantage of parallel computing architectures. Such approximation paradigm  gained a lot of popularity in the past years with the introduction of the so-called parareal method in the seminal series of work \cite{LionsMadayTurinici2001,MadayTuricini2002}. It consists in parallelizing the time integration of a dynamical system using an iterative process based on a predictor-corrector framework, but has never been used in the framework of collisional kinetic equation. Let us mention for the interested reader the very recent  book \cite{GanderLunet2024} that provides a comprehensive introduction to this topic. 

\smallskip

Let us now present the parareal approach for the simple ODE
\begin{equation*}
    \left\lbrace\begin{aligned}
    &\ddt u(t) = f(u),\quad t\in[0,T],\\
    &u(0) = u^0.
    \end{aligned}\right.
\end{equation*}
The goal of the parareal algorithm is to approximate the solution at some fixed discrete times $T^n$,
\begin{equation*}
    U^n \approx u(T^n),\quad n\in\{0,\dots,N_g\},\, N_g\in\N.
\end{equation*}
To do so, instead of using a single time integrator (or propagator), two solvers are coupled through an iterative process : the \emph{parareal iterations}. The first solver will be denoted by $\mathcal{G}(T^n,T^{n+1},U^n)$. It must be cheap to compute numerically, but yet convergent towards the solution to the ODE. It will provide some numerical \emph{guesses} of the behavior of the solution on the coarse grid. The second solver, that we will denote by $\mathcal{F}(T^n,T^{n+1},U^n)$, must be very accurate, and will be used to improve the aforementioned guesses.

While the propagator $\mathcal{G}$ is less costly, its drawback is naturally its precision. However, it will be corrected by the accurate, computationaly expensive, propagator $\mathcal{F}$. The algorithm then unfolds as follows:
\smallskip
\begin{enumerate}
    \item Divide the time domain in $N_g\in\N$ subintervals $[T^n,T^{n+1}]$, $n\in\{0,\dots,N_g-1\}$;
    \item Perform a first coarse guess: 
        \begin{equation*}
            U^{n+1,0} = \mathcal{G}(T^n,T^{n+1},U^{n,0}) = u^0;
        \end{equation*}
    \item Refine the guess through the parareal iterations:
    \begin{equation*}
        U^{n+1,k+1} = \mathcal{G}(T^n,T^{n+1},U^{n,k+1}) + \mathcal{F}(T^n,T^{n+1},U^{n,k}) - \mathcal{G}(T^n,T^{n+1},U^{n,k}),
    \end{equation*}
    where $k=1,2,\dots$ denotes the $k$\textsuperscript{th} parareal iteration.
\end{enumerate}
\smallskip
The advantage of this iterative algorithm is its ability to compute the expensive fine propagations in parallel. Therefore, as long as the number of parareal iterations remains small enough to obtain a given accuracy, one can expect a reduction of the computational cost of the time integration.

Parallel in time methods are now widely utilized and studied in various contexts. The convergence of the algorithm was investigated in \cite{GanderVandewalle2007,GanderHairer2008,DuarteMassotDescombe2011}, including variations of the original method. The uses are plenty, and we refer to \cite{BaudronLautardMadayMula2014} for neutron transport, \cite{FischerHechtMaday2005,SteinerRuprechtSpeckKrause2015} for the Navier-Stokes equations, \cite{GanderHairer2014} for Hamiltonian systems and \cite{SamaddarNewmanSanchez2010} for turbulence in plasmas. We also refer to a recent work on the acceleration of the method \cite{MadayMula2020, nguyen:tel-03950073}. It was furthermore observed that the parareal method suffers from poor convergence for hydrodynamic systems where convection dominates diffusion \cite{Gander2008,EghbalGerberAubanel2017,NielsenBrunnerHesthaven2018,Bal2005,FarhatChandesris2003}. In particular, in the case of Burgers' equation it was noticed that for a fixed discretization, the longer the integration domain, the slower was the convergence \cite{GanderHairer2008}.

In the context of multiscale systems, the parareal algorithm has been used in \cite{DuarteMassotDescombe2011,LegollLelievreSamaey2013,GrigoriHirstoagaNguyenSalomon2021,GrigoriSeverHirstoagaSalomon2023,SamaeySlawig2023,BossuytVandewalleSamaey2023} by leveraging reduced order models as coarse solvers. The idea is to use a simpler model, instead of a simpler time integrator, as a coarse propagator and a richer model for the fine propagations. These ideas are reminiscent to the ones introduced in the co-called Heterogeneous Multiscale Method framework from the series of works \cite{WeinanEngquist2003,Weinan2011PrinciplesOM,AbdulleWeinanEngquistVanden-Eijnden2012}.

\smallskip

Expanding on these advancements, our main goal in this work is to build a multiscale parareal method for kinetic equations to perform accurate simulations of the observable macroscopic moments. A fluid model, either given by a hydrodynamic or diffusive limit of the underlying kinetic equation, is used as a coarse propagator. The fine propagation is in turn achieved by the resolution of the kinetic model. To this aim, we resort to \emph{Asymptotic Preserving} (AP) schemes \cite{Klar1998,Jin1999} that remain stable for any value of the Knudsen number $\varepsilon$ in \eqref{para:VcolScaled}. While they provide an accurate approximation in the limit $\ep\rightarrow0$, their main drawback is that their cost remains the one of a kinetic scheme even in a fluid regime. We will demonstrate in our numerical experiments that our new method can solve this issue.

\medskip

\textbf{Plan of the paper.}
This work is organized as follows. In Section \ref{sec:MSParareal} we present the multiscale parareal algorithm for kinetic equations and discuss some of its properties. Sections \ref{para:sec:scheme} is dedicated to the implementation of the method and in particular to the presentation of the solvers used. We then discuss the parallelization of the method in Section \ref{para:sec:paraImp}. Finally, we thoroughly assess the accuracy and speedup properties of the method through numerical experiments in Section \ref{sec:NumResPara}.

\section{Multiscale parareal algorithm}\label{sec:MSParareal}
We introduce in this section the multiscale parareal algorithm. The formalism we shall employ follows the one introduced in \cite{LegollLelievreSamaey2013} and the references therein.

\subsection{Link between the scales}\label{para:sec:link}
The crucial point of a multiscale method is the link between the different scales of the model. In the context of kinetic equations a natural way to go from kinetic to fluid is to consider the projection of a distribution $f$ towards its moments $U=(\rho,u,\theta)^T$ defined in Definition \ref{para:defProj}.
\begin{Definition}\label{para:defProj}
    For a distribution function $f = f(t,x,v)\in L^1((1+v)^2\dd v)$, we define the projection $\mathcal{P}f(t,x)=U(t,x)$ as
    \begin{equation*}
    \mathcal{P}f = U = \begin{pmatrix}
        \int_{\R^3} f \dd v\\
        \\
        \frac{1}{\rho}\int_{\R^3} vf \dd v\\
        \\
        \frac{1}{3\rho}\int_{\R^3} |v-u|^2f \dd v\\
    \end{pmatrix}=\begin{pmatrix}
        \rho\\
        u\\
        \theta\\
    \end{pmatrix}.
    \end{equation*}
\end{Definition}
Note that in Definition \ref{para:defProj}, we start from a distribution that contains all the information and project it onto a manifold where the moments solve a system valid only for Maxwellians. On the opposite, to lift a macroscopic data $U(t,x)\in\R^5$ to a distribution $f(t,x,v)$ a first idea is to consider the validity of the fluid model used. In the case of Euler's equations, we therefore chose to reconstruct a Maxwellian whose moments are given by $U$.
\begin{Definition}\label{para:lift}
    For a macroscopic data $U(t,x)\in\R^5$, we define the lifting $\mathcal{L}U(t,x,v)=f(t,x,v)$ as
    \begin{equation*}
         \mathcal{L}U = f =\frac{\rho}{(2\pi\theta)^{3/2}}\exp\left(-\frac{|v-u|^2}{2\theta}\right).
    \end{equation*}
\end{Definition}

\begin{Remark}
Note that a finer way to handle the lifting and fluid propagations is to consider the Chapman-Enskog expansion of the distribution:
\begin{equation}\label{para:ChapmanEnskog}
    f^\ep=\M_{\rho^\ep,u^\ep\theta^\ep} + \sum_{l=1}^{\infty}\ep^l g^{(l)},
\end{equation}
where the perturbations $g^{(l)}$ can be explicitly computed and depend on $f$ only through its moments (some explicit computations for the BGK model we shall consider in this work can be fund in the book \cite{Struchtrup2005book}). 
\end{Remark}

We have already mentioned that in the hydrodynamic scaling, the moments of a Maxwel-lian distribution are solution to Euler's equations. By considering the first order perturbation in \eqref{para:ChapmanEnskog}, these moments can be shown to solve the Compressible Navier-Stokes equations \eqref{para:CNS} (CNS) that can be understood as a first order correction of Euler's equations. Consequently, a finer method would be to use a CNS solver as a coarse propagator. Note that even with this choice, the cost of the coarse integration remains much lower than a kinetic propagation. In addition, instead of only lifting the moments to a Maxwellian, one can consider adding the perturbations $g^{(l)}$. It allows for a finer reconstruction of the distribution and the consideration of far from equilibrium phenomena.
\begin{Definition}\label{para:lifthigher}
    With the same hypothesis of Definition \ref{para:lift}, we define the lifting
    \begin{equation*}
        \mathcal{L}^{(L)}U(t,x,v)=f(t,x,v)
    \end{equation*}
    of order $L$ by setting
    \begin{equation}
        f = \mathcal{L}^{(L)}U = \frac{\rho}{(2\pi\theta)^{3/2}}\exp\left(-\frac{|v-u|^2}{2\theta}\right) + \sum_{l=1}^{L}\ep^l g^{(l)}.
    \end{equation}
\end{Definition}

\subsection{Semi-discretization in time}
Let us now consider the time interval $[0,T]$ and divide it in $N_g\in\N$ uniform subintervals $[T^n,T^{n+1}]$. We are interested in the approximation of the moment vector $U$ at the fixed discrete times $T^n$. We denote by $U^{n,k}$ this approximation at a parareal iteration $k$. Furthermore, we also need to introduce a fine discretization that will be used by the fine propagator. Let $N_f\in\N$ denote the number of uniform fine subintervals $[t^n,t^{n+1}]$. Typically, $N_f$ shall be greater than $N_g$. We can then define the coarse (resp. fine) time steps and times:
\begin{equation*}
    T^n=n\Delta t_g,\quad \Delta t_g=\frac{T}{N_g},\quad  t^n=n\Delta t_f,\quad \Delta t_f=\frac{T}{N_f}.
\end{equation*}
Note that the output of the algorithm will therefore contain $N_g$ fixed snapshots at times $T^n$ of the moments. We briefly omit the definitions of the kinetic and fluid propagators $\mathcal{F}$ and $\mathcal{G}$ as they will be discussed below. However, it is important to mention that both the coarse and fine time steps are in practice the maximum time steps allowed within a propagator. The local time steps inside the solvers are naturally subject to a stability condition, and we define
\begin{equation*}
    \Delta t^{n,\mathcal{G}}_{\mathrm{loc}}=\min\left\lbrace\Delta t^{n,\mathcal{G}}_{\mathrm{stab}},\Delta t_g\right\rbrace,\qquad \Delta t^{n,\mathcal{F}}_{\mathrm{loc}}=\min\left\lbrace\Delta t^{n,\mathcal{F}}_{\mathrm{stab}},\Delta t_f\right\rbrace,
\end{equation*}
where $\Delta t^{n,\mathcal{G}}_{\mathrm{stab}}$ and $\Delta t^{n,\mathcal{F}}_{\mathrm{stab}}$ denote the time steps prescribed by the stability conditions of the numerical schemes used. The method can then be summarized in Algorithm \ref{para:algopara}.

\begin{algorithm}
    \caption{Multiscale kinetic parareal Algorithm}\label{para:algopara}
    \begin{algorithmic}[1]
    \Require $U^{0,0}$
    \For{$n = 1,\dots, N_g$} \Comment{First coarse guess}
        \State $U^{n,0}\gets \mathcal{G}\left(U^{n-1,0}\right)$ 
    \EndFor
    \vspace*{.2cm}
    \While{$k \leq K$ \textbf{or} error $\geq$  tol}\Comment{Parareal iterations}
    \vspace*{.2cm}
        \For{$n = 1,\dots, N_g$} \Comment{Compute the jumps in parallel}
            \State $\Delta^n = \mathcal{P}\mathcal{F}\mathcal{L}(U^{n-1,k-1}) - \mathcal{G}\left(U^{n-1,k-1}\right)$
        \EndFor
        \vspace*{.2cm}
        \For{$n = 1,\dots, N_g$} \Comment{Sequential correction}
            \State $U^{n,k+1} = \mathcal{G}\left(U^{n-1,k}\right) + \Delta^{n}$
        \EndFor
        \vspace*{.2cm}
        \State{Compute successive errors on the moments and $k\gets k+1$}
    \EndWhile
    \end{algorithmic}
\end{algorithm}

As a stopping criterion of the algorithm, we simply consider the error between two successive parareal iterations:
\begin{equation*}\label{para:conerrors}
    \mathrm{error} = \underset{n}{\max}{|U^{n,k+1}-U^{n,k}|},
\end{equation*}
or a fixed number of iterations $K\in\N$.

\smallskip 
An interesting property of Algorithm \ref{para:algopara} is that, at parareal iteration $k$, one can in fact start the loops (lines 5 and 8) at $k$ instead of $1$. Indeed, since one propagates a data from time $0$ that is fixed, it is unnecessary to recompute the first $1$ to $k$ iterations as they would not be modified. This can be seen through a simple ODE example by focusing on the first time intervals and parareal iterations: for $n=1$, to update $U^{1,k+1}$, one needs to compute
\begin{equation*}
    U^{1,k} = \mathcal{G}(T^0,T^{1},U^{0,1}) + \mathcal{F}(T^0,T^{1},U^{0,0}) - \mathcal{G}(T^0,T^{1},U^{0,0}).
\end{equation*}
Since $U^{0,k}=U^{0}$ for all $k$, the two coarsely propagated data cancel each other, and we obtain that the first parareal iteration automatically corrects the initial guess towards the fine propagation:
\begin{equation*}
    U^{1,k} = \mathcal{F}(T^0,T^{1},U^{0,0}),\quad\forall k\geq1.
\end{equation*}
Consequently, when computing the next parareal iteration $k=2$, $U^{1,k}$ is now fixed and one can make the same observation as before but now for $n=2$: 
\begin{equation*}
    U^{2,k} = \mathcal{F}(T^1,T^{2},U^{1,k}) = \mathcal{F}(T^1,T^{2},\mathcal{F}(T^0,T^{1},U^{0,0})),\quad\forall k\geq2.
\end{equation*}
By induction, at parareal iteration $k$, the $1$ to $k$ first times are therefore fixed, and don't need to be recomputed. A consequence is that after $N_g$ parareal iterations, the outcome of Algorithm \ref{para:algopara} is the same as a fully fine simulation.

Building upon these observations, an optimized version of this algorithm is given by Algorithm~\ref{para:algoparaBis}.

\begin{algorithm}
    \caption{Optimized multiscale kinetic parareal Algorithm}\label{para:algoparaBis}
    \begin{algorithmic}[1]
    \Require $U^{0,0}$
    \For{$n = 1,\dots, N_g$} \Comment{First coarse guess}
        \State $U^{n,0}\gets \mathcal{G}\left(U^{n-1,0}\right)$ 
    \EndFor
    \vspace*{.2cm}
    \While{$k \leq K$ \textbf{or} error $\geq$  tol}\Comment{Parareal iterations}
    \vspace*{.2cm}
        \For{$n = k,\dots, N_g$} \Comment{Compute the jumps in parallel}
            \State $\Delta^n = \mathcal{P}\mathcal{F}\mathcal{L}(U^{n-1,k-1}) - \mathcal{G}\left(U^{n-1,k-1}\right)$
        \EndFor
        \vspace*{.2cm}
        \For{$n = k,\dots, N_g$} \Comment{Sequential correction}
            \State $U^{n,k+1} = \mathcal{G}\left(U^{n-1,k}\right) + \Delta^{n}$
        \EndFor
        \vspace*{.2cm}
        \State{Compute successive error on the moments and $k\gets k+1$}
    \EndWhile
    \end{algorithmic}
\end{algorithm}

From an optimization point of view, note also that steps 5-7 of both Algorithms \ref{para:algopara} and \ref{para:algoparaBis} contain all the expensive computations. More particularly we chose to do the lifting and the projection alongside the kinetic propagation to better distribute the workload. Indeed, since the iterations of this loop are independant from one another, which is the core of the parareal method, it can be parallelized to reduce the computation time.

\section{Numerical schemes}\label{para:sec:scheme}
In this section, we shall detail our choice of fine and coarse propagators as well as the discrete versions of the lifting and projection operators. 

From now on, let us place ourselves in the $d_x=1$ and $d_v=3$ setting that we shall use for our numerical experiments. The $(x,v)$--phase space  is discretized in a finite volume fashion using control volumes $K_{ij}=\mathcal{X}_i\otimes\mathcal{V}_j$, for $i\in\{1,\dots,N_x\}$ and $j=(j_x,j_y,j_z)\in\{1,\dots,N_v\}^3$. The discrete unknowns are defined as
\begin{equation*}
    f_{ij}^n \approx \frac{1}{\Delta x_i\Delta v_j}\int_{K_{ij}} f(T^n,x,v)\dd x \dd v,\quad U_i^n \approx \frac{1}{\Delta x_i}\int_{\mathcal{X}_i} U(T^n,x)\dd x,
\end{equation*}
with $\Delta x_i$ the volume of $\mathcal{X}_i$ and $\Delta v_j$ the volume of the cube $\mathcal{V}_J$.

\subsection{Kinetic schemes}
The main difficulty to approximate scaled kinetic equations is the stability of the time integration with respect to the small parameters. In both hydrodynamic and diffusive scalings, AP schemes must therefore be considered to ensure tractable simulations for any value of the Knudsen number.

\medskip

\textbf{Hydrodynamic scaling.} In the hydrodynamic scaling, several AP schemes have been developed, and most of them rely on an implicit-explicit (IMEX) type time discretization \cite{PareschiRusso2005}. In this work, we consider a very simple IMEX discretization where the stiff relaxation term is implicit. The semi-discretized scheme can be written as:
\begin{equation}
    \frac{f^{n+1} - f^n}{\Delta t} + \left[\, v \cdot \nabla_x f^n + E\cdot \nabla_v f^n \,\right] = \frac{\tau^{n+1}}{\varepsilon} \left(\M_{\rho^{n+1},u^{n+1},\theta^{n+1}} - f^{n+1}\right).
\end{equation}
Thanks to the conservation properties of the equation, the moments $\rho$, $u$ and $\theta$ at time $n+1$ used to build the Maxwellian $\M^{n+1}$ can be obtained explicitly, yielding an efficient, fully explicit in time procedure \cite{CoronPerthame1991, FilbetJin2010}. More advanced methods such as a \textit{micro-macro} approach could also be considered here \cite{LemouMieussens2008}. For the phase space discretization, the numerical fluxes are of upwind type for the transport in position and of Rusanov type for the transport in velocity \cite{Leveque2002,Toro2009}. The fully discretized kinetic system now reads
\begin{equation}
    \frac{f_{ij}^{n+1} - f_{ij}^n}{\Delta t} + \left[\, \mathcal{F}_{\iph,j}^n-\mathcal{F}_{\imh,j}^n + \mathcal{G}_\jph^n-\mathcal{G}_\jmh^n \,\right] = \frac{\tau_i^{n+1}}{\varepsilon} \left(\M_{\rho_i^{n+1},u_i^{n+1},\theta_i^{n+1}} - f_{ij}^{n+1}\right),
\end{equation}
with the numerical fluxes
\begin{equation}
    \begin{aligned}
        &\mathcal{F}_{\iph,j}^n= v^+ f_{ij}^n + v^- f_{i+1,j}^n,
        &\mathcal{G}_{i,\jph}^n= \frac{E_i}{2}(f_{ij}^n + f_{i,j+1}^n) + \frac{\underset{i}{E_{max}}}{2}(f_{i,j+1}^n - f_{i,j+1}^n),
    \end{aligned}
\end{equation}
where $a^+=\max\{a,0\}$, $a^-=-\min\{a,0\}$ and $E_{max} = \max{E_i}$.

It is widely known that such an approach yields a stable, yet first order discretization of the collisional kinetic equation \eqref{para:VcolScaled}.
A natural improvement to this fine propagator would be to use fluxes limiters \cite{van1979towards}.

\medskip

\textbf{Diffusive scaling.} In the diffusive scaling, one could adopt the micro-macro approach introduced in \cite{LemouMieussens2008,CrouseillesLemou2011} together with an exponential time integrator for the microscopic equation \cite{Lemou2010,Laidin2023}. This approach will be explored in a future work.

\subsection{Fluid schemes}
Coarse propagators aim at giving a decent guess of the behavior of the moments at a low numerical cost. We distinguish two types of schemes, depending on the scaling considered at the kinetic level.

\medskip

\textbf{Systems of conservation laws.} In the framework of multiscale kinetic equations, systems of conservation laws typically arise as hydrodynamic limits of kinetic equations. They can generally be written as:
\begin{equation}
    \dt U + \nabla_x\cdot F(U) = S(U),
\end{equation}
where $F$ is the flux function and $S$ a potential source term  \cite{Leveque2002,Toro2009}. As a coarse propagator, we rely on a simple finite volume discretization with Rusanov-type fluxes. 

The Euler equations in this setting can be expressed as
                \begin{equation*}
                    \frac{\partial V}{\partial t} + \frac{\partial {h}(V)}{\partial x}  =0,
                \end{equation*}
                where
                \begin{equation*}
                    V=
                    \begin{bmatrix}
                    \rho \\
                    \rho u_x\\
                    \rho u_y\\
                    \rho u_z \\
                    \mathcal E\\
                    \end{bmatrix}, \ \ 
                     {h}(V)=
                    \begin{bmatrix}
                    \rho u_x \\
                    \rho u_x^2 + \rho \theta\\
                    \rho u_x u_y\\
                    \rho u_x u_z\\
                    u_x (\mathcal E+\rho \theta)\\
                    \end{bmatrix}
                \end{equation*}
                One can use the following explicit formula  to compute the coarse integrator
                \begin{equation}
                    V^{n+1}_{i,j}=V^{n}_{i,j} -         \left[H_{i+\frac{1}{2}}^n - H_{i-\frac{1}{2}}^n\right],
                \end{equation}
                where the Rusanov  flux  $H$ is given by
                \begin{align*}
                    H_{i+\frac{1}{2}}^n &=
                    \frac{1}{2}\left({h}(V_{i+1}^n) - {h}(V_{i}^n)\right) - \frac{1}{2} S_{i+\frac{1}{2}}^n \left(U_{i+1}^n - U_{i}^n\right),
                \end{align*}
                with
                \begin{align*}
                    S_{i+\frac{1}{2}}^n &=
                    \max\left(|u^n_{x, i+1}|+c_{i+1}, |u^n_{x, i}|+c^n_{i,j}\right), 
                \end{align*}
                and $c_{i}^n=\sqrt{\theta^n_{i}}$ is the sound speed.

\medskip

\textbf{Drift-diffusion equation.} Explicit discretizations of diffusive systems usually suffer from a parabolic stability condition. However, since the cost of such a scheme is essentially negligible compared to a full kinetic simulation, we do not consider an implicit-in-time procedure. The fluid system is therefore approximated using a central finite difference scheme that is nothing but the asymptotic limit of the micro-macro scheme mentioned above \cite{Lemou2010,Laidin2023}.

\begin{Remark}
    As is the case for the fine propagator, it is common to consider the use of slope of flux limiters when solving systems of conservation laws. Nevertheless, we noticed in practice that additional spatial accuracy on the predictor does not accelerate the convergence of the parareal algorithm. On the contrary, a slower convergence was observed. 
    Our interpretation is that the numerical viscosity is beneficial. This corroborates with the fact the parareal algorithm is known to work best when the solution is regular, as was already mentioned in the series of works \cite{Gander2008,EghbalGerberAubanel2017,NielsenBrunnerHesthaven2018,Bal2005,FarhatChandesris2003}.
\end{Remark}

\subsection{Discrete lifting and projection operators}
Let us finish the complete description of our numerical propagators by describing the discrete lifting and projection operators.

\medskip
\textbf{Discrete lifting operator.}
A natural way to define the discrete lifting operator is to compute the pointwise discrete Maxwellian associated to the moments $U_i=(\rho_i,u_i,\theta_i)^T$, $i\in\{1,\dots,N_x\}$, namely
\begin{equation}
    (\mathcal{L}U)_{ij} = \frac{\rho_i}{\left(2\pi\theta_i\right)^{3/2}}\exp\left(-\frac{|v_j-u_i|^2}{2\theta_i}\right),\quad j=(j_x,j_y,j_z)\in\{1,\dots,N_v\}^3.
\end{equation}
We could also define higher order lifting as defined in the continuous case \eqref{para:lifthigher} and proceed in the same way using pointwise computations and finite difference approximations of the derivatives. 

\medskip
\textbf{Discrete projection operator.}
To project a given discrete distribution $f_{ij}$ defined on a discrete phase space towards its discrete moments $(\rho_i,u_i,\theta_i)^T$, we consider a simple first order quadrature:
\begin{equation*}
    \rho_i=\sum_{j} f_{ij}\Delta v_j,\quad (\rho u)_i=\sum_{j} v_j f_{ij}\Delta v_j,\quad3(\rho\theta)_i=\sum_{j} |v_j-u_i|^2f_{ij}\Delta v_j.
\end{equation*}
Note that more accurate quadratures could be considered depending on the velocity discretization and a desired accuracy.

\section{Parallelization of the method}\label{para:sec:paraImp}
In order to achieve a significant gain in computational time, one needs to rely on an efficient parallelization. The method is currently implemented using an OpenMP paradigm to validate its concept. Its real performance should stand out when properly deployed on an heterogeneous supercomputer along with an MPI parallelization over several threads. The parallelization of this type of method on distributed memory architecture is known to be challenging and such an implementation is briefly discussed in Appendix \ref{AppendixMPIimpl}.

\medskip

\textbf{Theoretical speed up.}
In an ideal setting, \textit{i.e.} by omitting the cost of communications and synchronizations between cores/threads, one can derive an estimate on the number of parareal iterations to obtain a speedup with PinT-types methods. 

Let us denote by $T_{\mathrm{Kin}}$ (resp. $T_{\mathrm{Fluid}}$) the maximum time for the fine, kinetic (resp. coarse, fluid) solver to evolve an initial data from time $T^n$ to time $T^{n+1}$. Note that we take the maximum over all propagations because, on each sub time-interval, the stability condition may vary, impacting the optimal local time step. In addition, one needs to take into account the cost of the projection and lifting operators since they are fully non-local and can therefore be significantly expensive. We shall denote by $T_{\mathrm{Lift}}$ and $T_{\mathrm{Proj}}$ these costs. The ideal numerical cost of Algorithm \ref{para:algopara} on $N_p$ threads with $N_g$ sub-time intervals is then given by the formula:
\begin{equation*}
    T_{\mathrm{Parareal}} = T_{\mathrm{Fluid}} + N_g k \left(\frac{T_{\mathrm{Lift}}+T_{\mathrm{Proj}}+T_{\mathrm{Kin}}+T_{\mathrm{Fluid}}}{N_p} + T_{\mathrm{Fluid}}\right).
\end{equation*}
Consequently, the ideal number of parareal iterations to be less costly than a fully kinetic simulation should satisfy:
\begin{equation*}
    T_{\mathrm{Parareal}}\leq N_gT_{\mathrm{Kin}}
\end{equation*}
or stated differently,
\begin{equation*}
    k_{\mathrm{opt}} \leq \left\lceil\frac{N_gT_{\mathrm{Kin}}-T_{\mathrm{Fluid}}}{N_g\left(\frac{T_{\mathrm{Lift}}+T_{\mathrm{Proj}}+T_{\mathrm{Kin}}+T_{\mathrm{Fluid}}}{N_p}+T_{\mathrm{Fluid}}\right)}\right\rceil.
\end{equation*}

\medskip
\textbf{Practical implementation.}
The performance of the method and scalability of the current implementation for this paper were made in Fortran90 (gfortran 9.4.0) with an OpenMP parallelization that will be discussed in the Section \ref{sec:NumResPara}. These tests were ran on the architecture detailed in Table \ref{para:architecture}.

\begin{table}
    \centering
    \begin{tabular}{c|c}
        \# CPUs &  Intel(R) Xeon(R) Platinum 8268 CPU @ 2.90GHz x2\\
        \hline
        \# cores & 48 \\
        \hline
        RAM & 750 GB \\
    \end{tabular}
    \caption{Computing architecture.}
    \label{para:architecture}
\end{table}
In a shared memory paradigm, the parallelization of the parallel loop in algorithm \ref{para:algoparaBis} is achieved using a dynamic scheduling as the number of iteration is not fixed.

\section{Numerical results}\label{sec:NumResPara}
In this section, we present some numerical results for simulations of the 1D/3D Vlasov-BGK equation \eqref{para:VcolScaled} in the hydrodynamic scaling. We consider Algorithm \ref{para:algoparaBis} and assess its accuracy as well as its performance through several experiments. In the following, unless specified otherwise, the phase space is uniformly discretized using $200\times32\times32\times32$ cells, namely $200$ points in position and $32$ in each velocity direction. Note that the method was only implemented in the hydrodynamic scaling with an Euler solver at the fluid level and the lifting does not consider any perturbations and then any higher order models. Extensions to the cases mentioned in Section \ref{para:sec:scheme} will be implemented in future works.

\subsection{Test 1: Sod shock tube}
We first consider a Riemann problem where the initial data is given by:
\begin{equation*}
    f_0(x,v) = \M_{\rho(x),u(x),\theta(x)}(v),\quad x\in[0,2], \quad v\in[-8,8],
\end{equation*}
where
\begin{equation*}
    (\rho(x),u(x),\theta(x))=\left\lbrace\begin{aligned}
        &(1,0,0,0,1) \quad&\text{if }x<1,\\
        &(0.125,0,0,0,0.8) &\text{if } x\geq1.
    \end{aligned}\right.
\end{equation*}
The exterior force is set to $0$, and the regime is rarefied: $\ep=10^{-2}$ at the kinetic level. Boundary conditions are of the absorbing type. We set a simulation interval $[0,T]$ with $T=0.5$ and discretize it uniformly using $N_g=200$ points for the coarse grid and $N_f=800$ points for the fine grid.  As a stopping criterion, we set a maximum of $K=80$ iterations or when then consecutive error \eqref{para:conerrors} is smaller than $10^{-8}$. Note that this threshold is much smaller than both the fine and coarse time increments.

\begin{figure}
    \centering
    \includegraphics[width=.99\linewidth]{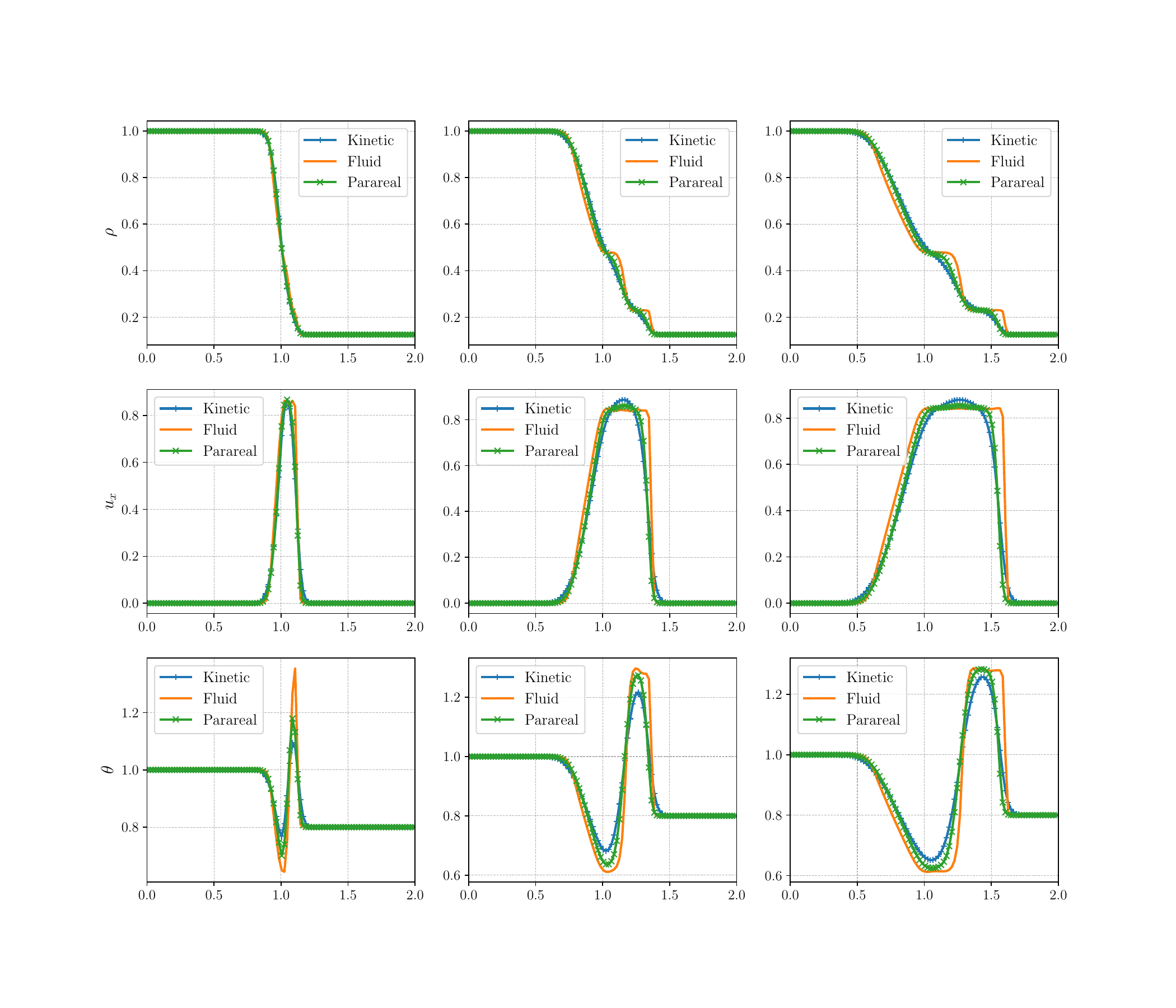}
    \caption{\textbf{Test 1 - Sod shock tube}, $\ep=10^{-2}$: Snapshots of the density (Top), $x$ mean velocity (Middle) and Temperature (Bottom) at times $T^n=0.05$ (Left), $0.15$ (Middle) and $0.25$ (Right).}
    \label{para:fig:SODSnap0.01}
\end{figure}

We present in Figure \ref{para:fig:SODSnap0.01} snapshots of the moments obtained with Algorithm \ref{para:algoparaBis} as well as the ones obtained through fully kinetic and fully fluid simulations. We can observe a very good agreement between the kinetic and parareal densities, even if the results of the fully fluid equation (that we recall is used as coarse propagator) are far from exact. Note that we also notice that higher order moments are less well captured, and this observation also holds as we look for larger times. Nevertheless, considering we used the simplest fluid model and reconstruction, the results are promising.

\medskip

\begin{figure}
    \centering
    \includegraphics[width=.60\linewidth]{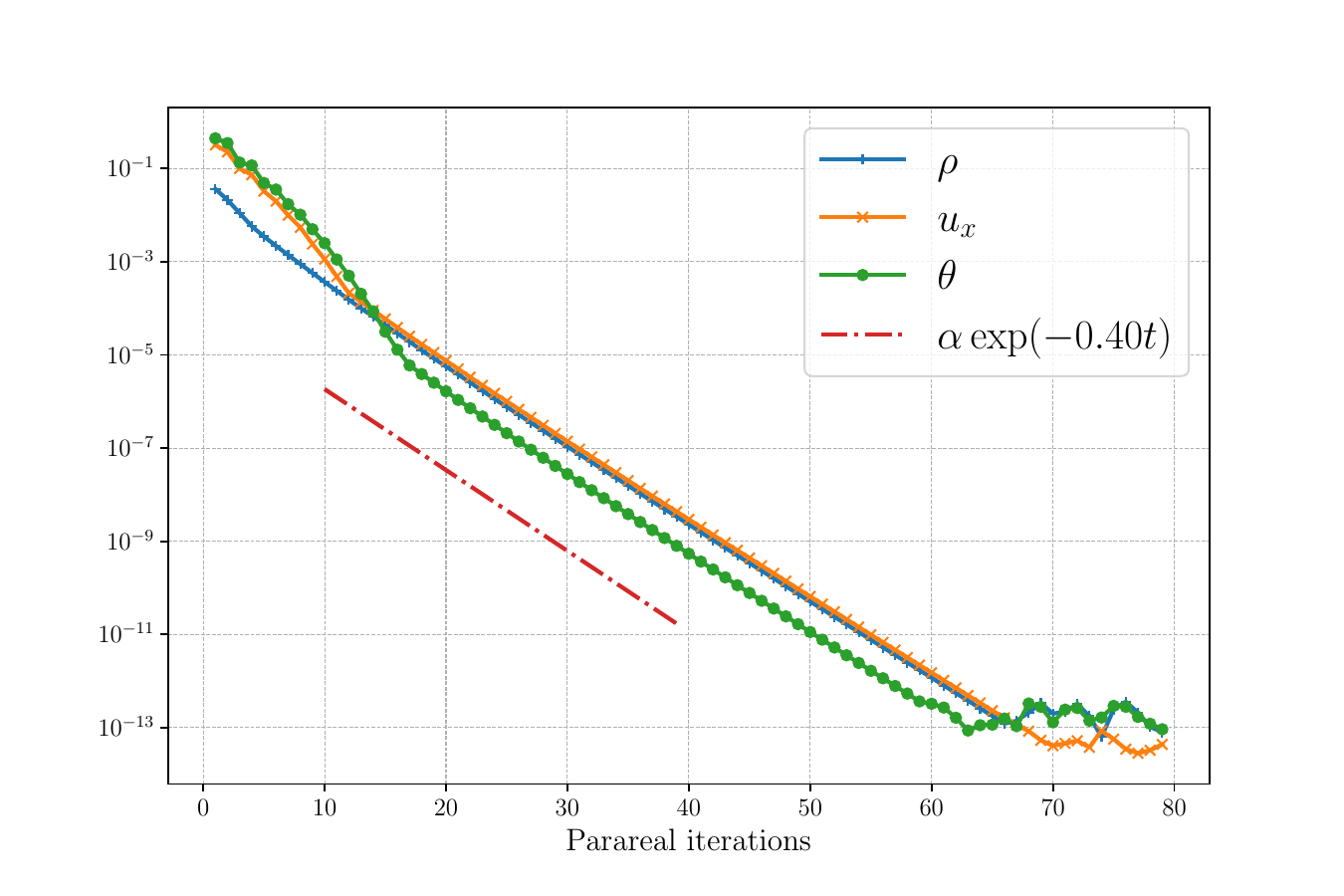}
    \caption{\textbf{Test 1 - Sod shock tube}, $\ep=10^{-2}$: Convergence of the successive errors.}
    \label{para:fig:SODconv0.01}
\end{figure}

\textbf{Convergence.} Figure \ref{para:fig:SODconv0.01} illustrates the convergence of the algorithm, by plotting in loglog scale the successive errors \eqref{para:conerrors} at each parareal iteration. A first observation is that the algorithm indeed converges well. In particular, the error decreases exponentially fast with the number of parareal iterations, which is encouraging in the sense that one may need only a few iterations to reach a desired accuracy.

\subsection{Test 2: Blast waves}
We now consider a second Riemann problem consisting of two hot, energetic shock waves converging asymmetricaly towards a same point, namely a blast wave-type collision. The initial data is at a local equilibrium in velocity.
\begin{equation*}
    f_0(x,v) = \M_{\rho(x),u(x),\theta(x)}(v),\quad x\in[0,2], \quad v\in[-8,8],
\end{equation*}
with the moments defined by:
\begin{equation*}
    (\rho(x),u(x),\theta(x))=\left\lbrace\begin{aligned}
        &(1,1,0,0,2) \quad&\text{if }x<0.4,\\
        &(1,0,0,0,0.25) &\text{if } 0.4\leq x<1.6,\\
        &(1,-1,0,0,2) &\text{if } x\geq1.6.
    \end{aligned}\right.
\end{equation*}

\begin{figure}
    \centering
    \includegraphics[width=.99\linewidth]{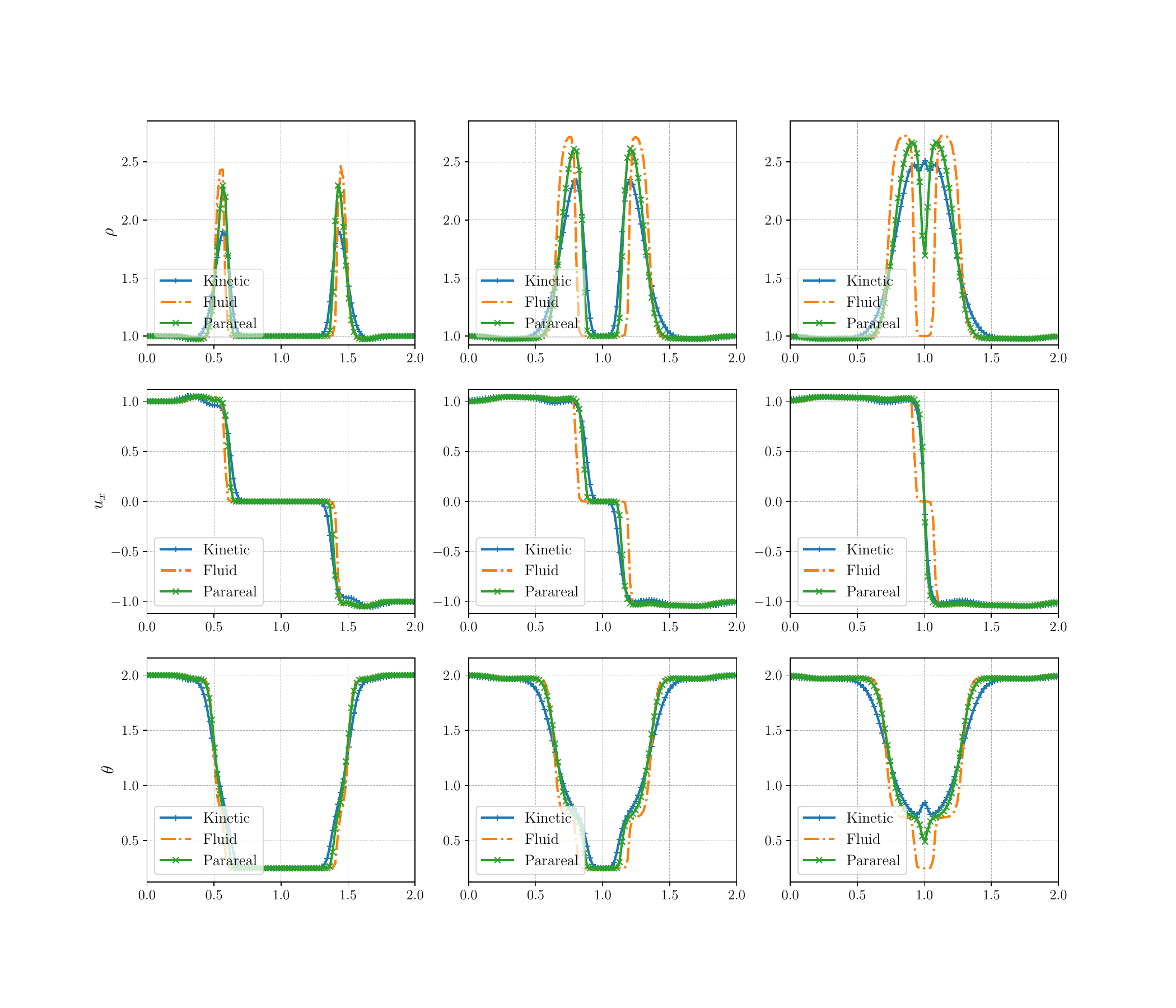}
    \caption{\textbf{Test 2 - Blast waves}, $\ep=10^{-2}$: Snapshots of the density (Top), $x$ mean velocity (Middle) and Temperature (Bottom) at times $T^n=0.1$ (Left), $0.23$ (Middle) and $0.3$ (Right).}
    \label{para:fig:BlastSnap0.01}
\end{figure}

We consider the same setting as in Test 1, and we fix the number of parareal iterations to $K=10$. We present in Figure \ref{para:fig:BlastSnap0.01} snapshots of the moments. Note that with only 10 iterations, the consecutive error is only of the order of $10^{-3}$. We observe that the solution has not yet converged towards the kinetic one, and this is even more striking for large time. Nevertheless, this behavior is expected as the parareal algorithm first corrects for early times and then propagates this correction.

\begin{figure}
    \centering
    \includegraphics[width=.99\linewidth]{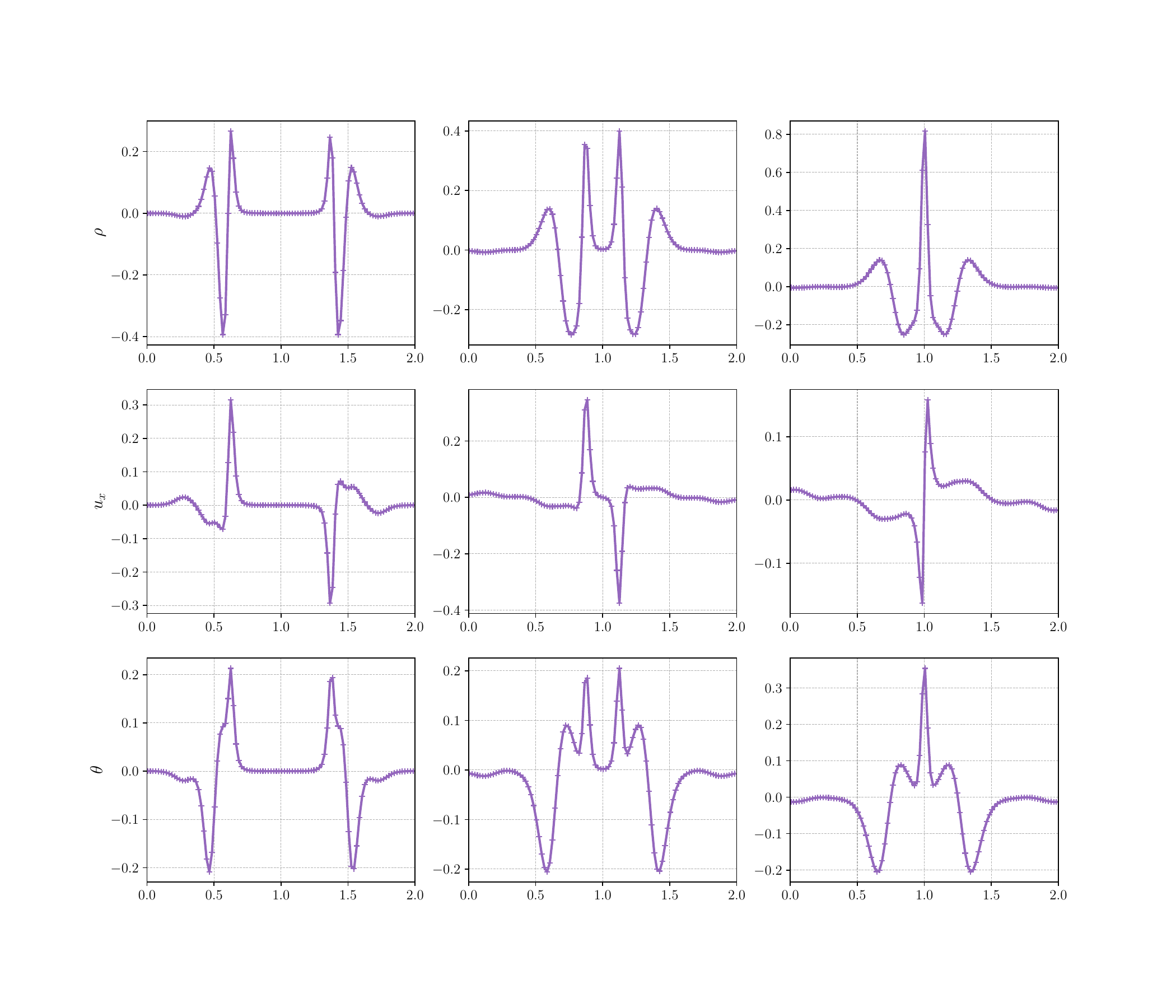}
    \caption{\textbf{Test 2 - Blast waves}, $\ep=10^{-2}$: Snapshots of the pointwise difference on the density (Top), $x$ mean velocity (Middle) and Temperature (Bottom) at times $T^n=0.1$ (Left), $0.23$ (Middle) and $0.3$ (Right).}
    \label{para:fig:BlastSnapError0.01}
\end{figure}

In Figure \ref{para:fig:BlastSnapError0.01}, we investigate the main source of error by plotting the point-wise difference between the kinetic and parareal moments. According to Figure \ref{para:fig:BlastSnap0.01} we observe that the deviation between kinetic and parareal solvers is essentially localized near shocks, where the regularity of the solution is low. Again, such behavior is not surprising as it was already observed in previous works \cite{Gander2008,EghbalGerberAubanel2017,NielsenBrunnerHesthaven2018,Bal2005,FarhatChandesris2003}.

\begin{figure}
    \centering
    \includegraphics[width=.60\linewidth]{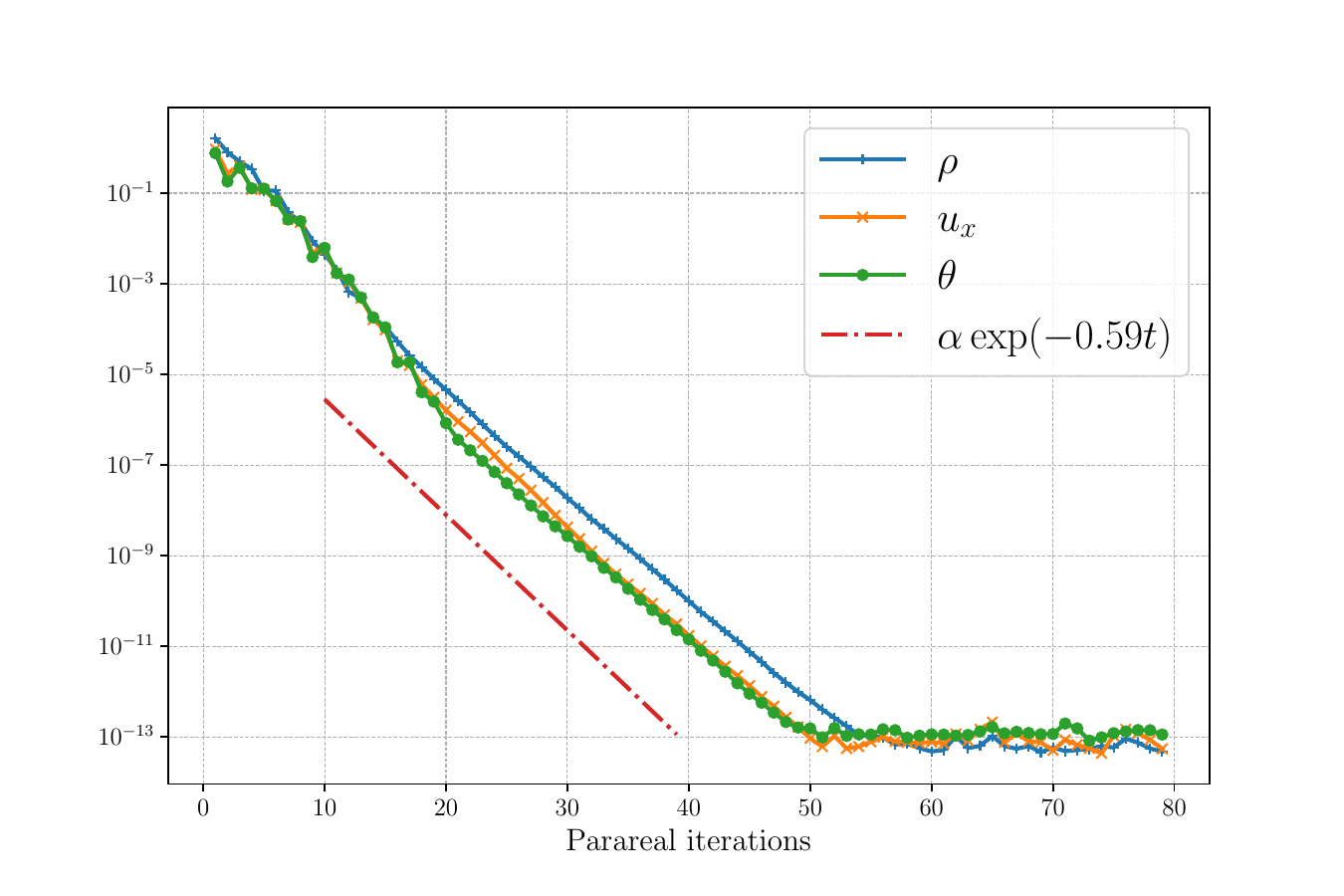}
    \caption{\textbf{Test 2 - Blast waves}, $\ep=10^{-2}$: Convergence of the successive errors.}
    \label{para:fig:Blastconv0.01}
\end{figure}

We also report in Figure \ref{para:fig:Blastconv0.01} the convergence of the algorithm where we observe the same exponetial decay with respect of the number of parareal iteration, as was noticed in Test 1.

\subsection{Test 3: Exterior forces}
\begin{figure}
    \centering
    \includegraphics[width=.50\linewidth]{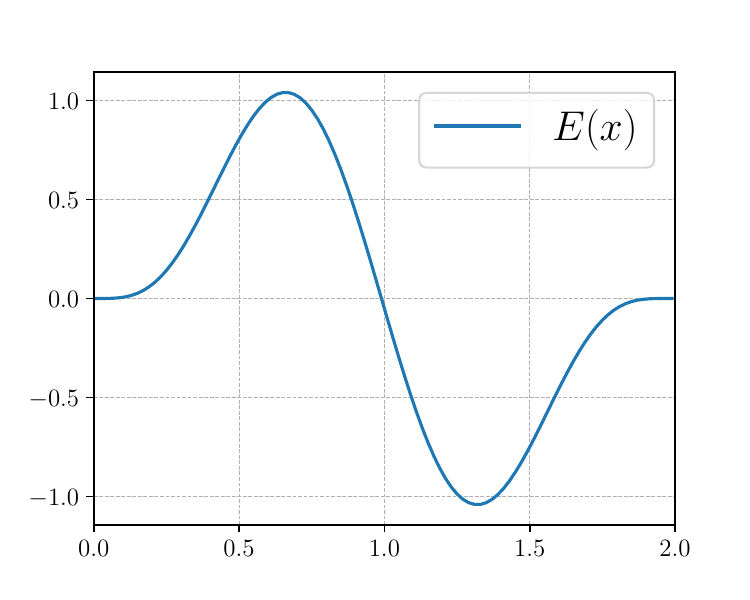}
    \caption{\textbf{Test 3 - Exterior forces:} Plot of the external force $E(x)$.}
    \label{para:fig:PlotExtForce}
\end{figure}

We conclude this numerical investigation of the multiscale PinT method by considering the case of a smooth initial data along with a non-zero exterior force. Moreover, we consider this time an initial data far from the Maxwellian equilibrium given by two opposite beams with constant densities:
\begin{equation*}
    f_0(x,v) = \M_{\rho_1(x),u_1(x),\theta_1(x)}(v) + \M_{\rho_2(x),u_2(x),\theta_2(x)}(v),\quad x\in[0,2], \quad v\in[-8,8],
\end{equation*}
where the two sets of moments are given by:
\begin{equation*}
    (\rho_1(x),u_1(x),\theta_1(x))=(1,1,0,0,1)
\end{equation*}
and 
\begin{equation*}
    (\rho_2(x),u_2(x),\theta_2(x))=(1,-1,0,0,1).
\end{equation*}
We set the exterior forces as
\begin{equation*}
    E(x) = -5x^4(x-2)^4(x-1).
\end{equation*}
This exterior force is illustrated in Figure \ref{para:fig:PlotExtForce} and will concentrate the mass towards the center of the domain $x=1$.
A natural follow up to this test will be to consider a coupling with Poisson's equation.

The discretization is chosen as $100\times256\times16\times16$ cells, namely $100$ points in position, $256$ points in the $v_x$ directions and $32$ in the $v_y$ and $v_z$ directions. The specific refinement in $v_x$ is only meant to reduce the consistency error of our kinetic solver  as there now is transport in this particular direction. Note that the number of cells in the other directions is unchanged compared to the previous test cases. In addition, we assume periodic boundary conditions and a Knudsen number of $10^{-5}$. 

\begin{figure}
    \centering
    \includegraphics[width=.99\linewidth]{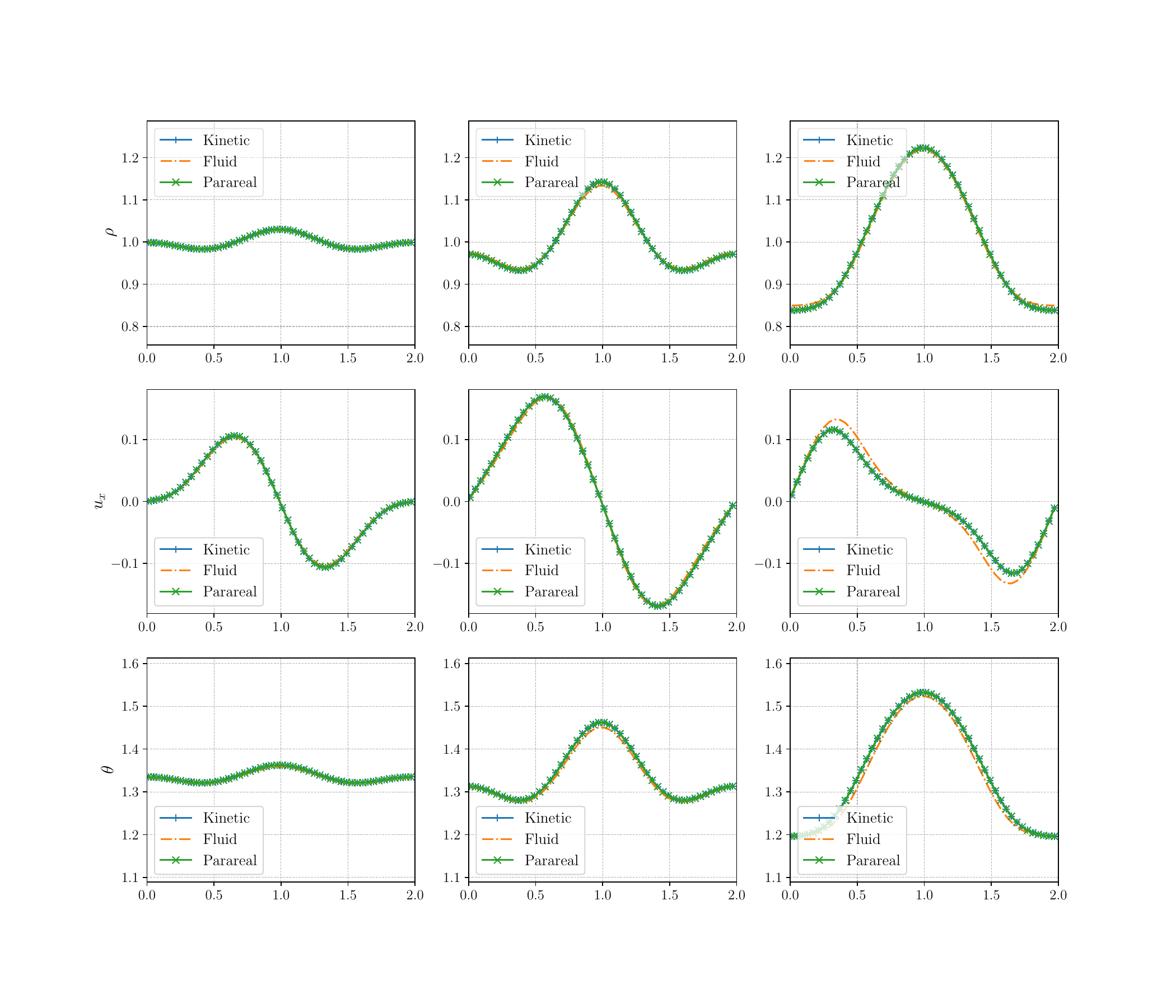}
    \caption{\textbf{Test 3 - Exterior forces}, $\ep=10^{-5}$: Snapshots of the density (Top), $x$ mean velocity (Middle) and Temperature (Bottom) at times $T^n=0.1$ (Left), $0.25$ (Middle) and $0.5$ (Right).}
    \label{para:fig:SmoothSnap1e-5}
\end{figure}

The resulting moments are presented in Figure \ref{para:fig:SmoothSnap1e-5}. We can observe that the mass indeed concentrates towards the center of the domain because of the exterior force. In particular, the solution being smooth, we obtain a very good agreement between kinetic and parareal moments. In addition, we also illustrate the convergence of the algorithm in the next Paragraph along with the performance of the method. 

\medskip

\textbf{Performance.} We conclude this test by first investigating the parallel efficiency of our implementation. The scalability is assessed by measuring the time to perform the first, most expensive parareal iteration. The results obtained on the architecture presented in Table \ref{para:architecture} are illustrated in Table \ref{tab:Scalability}. The current implementation scales as around $0.737$ with the number of threads which is very good considering the parallelization method used. A refinement of this implementation in MPI is discussed in Appendix \ref{AppendixMPIimpl}.

\begin{table}
    \centering
    \begin{tabular}{c|c|c|c|c|c|c}
        \# Threads &  1 & 2 & 4 & 8 & 16 & 32\\
        \hline
        Runtime (s) & 278.7 & 153.5 & 84.3 & 48.9 & 30.6 & 22.9
    \end{tabular}
    \caption{\textbf{Test 3 - Exterior forces}, $\varepsilon=10^{-5}$:}
    \label{tab:Scalability}
\end{table}

Finally, we report in Table \ref{tab:Speedup} the successive error, the runtimes and the speedup obtained with our new method. Firstly, we observe a significant reduction of the computational time with a factor that naturally depends on the the number of parareal iterations. For a reasonable error of $10^{-6}$, the observed speedup reaches $4.52$ between the full, direct kinetic integration, and the new multiscale PinT one. We believe that this is  particularly significant, specially when one considers our simplified $d_x=1$ setting.

\begin{table}
    \centering
    \begin{tabular}{c|c|c|c|c}
      Parareal Iterations & Final error & Runtime (s) & Kinetic Runtime & Speedup \\
      \hline
      3  & 1.387E-03 & 71.7 & 672.8 & 9.38 \\
      \hline
      6  & 7.782E-06 & 153.0 & 672.8 & 4.52 \\
      \hline
      9 & 2.421E-08 & 218.4 & 672.8 & 3.08
    \end{tabular}
    \caption{\textbf{Test 3 - Exterior forces},  $\varepsilon=10^{-5}$: Performance of Algorithm \ref{para:algoparaBis} on 32 threads.}
    \label{tab:Speedup}
\end{table}

\section{Conclusion}
We presented in this work the first multiscale, parareal in time numerical solver for computing efficiently solution to collisional, Vlasov-type kinetic equation. We introduced our method, analysed some of its numerical properties, proposed a first implementation in a shared memory framework and conducted a thorough testing of this new method. We exhibited some very serious speedup (up to 10 times faster than a full kinetic run) using a regular, not state-of-the-art computer. We finally presented a possible distributed memory version of our algorithm that we hope to test in bigger computational clusters in future works.

In addition, it is worth mentioning that we expect even larger computational gains in more complex cases such as higher dimensions in positions and Boltzmann-type collision operators. Moreover, as it was for example observed in \cite{MadayMula2020}, a speedup can also be more significant for large time simulations, which have not been considered in this work. 

On a final note, more sophisticated versions of the parareal algorithm have been recently constructed with improvement of performance in mind and we refer to \cite{MadayMula2020,nguyen:tel-03950073} for more details that we also wish to consider in future works.

\section*{Acknowledgement}
The authors would like to warmly thank M. Bessemoulin-Chattard for the careful reading and feedback on early versions of the manuscript that happened at the end of T. Laidin's PhD. 
T. Laidin has received funding from the European Research Council (ERC) under the European Union’s Horizon 2020 research and innovation program (grant agreement No 865711).
T. Rey received funding from the European Union's Horizon Europe research and innovation program under the Marie Sklodowska-Curie Doctoral Network DataHyKing (Grant No. 101072546), and is supported by the French government, through the UniCA$_{JEDI}$ Investments in the Future project managed by the National Research Agency (ANR) with the reference number ANR-15-IDEX-01.

\appendix

\section{MPI implementation of the multiscale parareal method}\label{AppendixMPIimpl}
This appendix is dedicated to proposing an MPI implementation of Algorithm \ref{para:algoparaBis}. We present the procedure to be implemented and discuss the challenges it poses.

\medskip

\textbf{Workload distribution.} Let us consider the parallelization of Algorithm \ref{para:algoparaBis} over $N_p\in\N$ processors with distributed memory. We recall that we aim at computing the expensive kinetic propagations in parallel to reduce the simulation time. 

A first observation is that after each parareal iteration, the actual work to be done diminishes. Therefore, one must adjust the load distribution at each iteration. Let us now consider a discretization of the position variable using $N_x$ points and $N_g$ points for the coarse time discretization. From a practical point of view, the memory allocation of the jumps arrays ($\Delta^n$ in Algorithm \ref{para:algoparaBis}) should be done locally on each processor, and only once, to avoid multiple allocation/deallocation overheads. Indeed, since they depend on the size of the spatial discretization, their memory footprint is potentially large. Therefore, the master processor shall be allocated with the full time discretization which is of size $5N_xN_g$, where $5$ corresponds to the number of moments. The remaining processors are allocated their maximum workload, namely their workload at the first parareal iteration. This quantity, denoted by $\texttt{chunk}$ in the following, can be defined using integer arithmetic to evenly distribute the time iterations. 

The workload at parareal iteration $k$, the starting and ending time indices for each processor are then computed as 
\begin{equation*}
    \begin{aligned}
        &\text{\texttt{work}} = N_g-k;\\
        &\text{\texttt{chunk}} = \left\lfloor\text{\texttt{work}} / N_p\right\rfloor;\\
        &\text{\texttt{remainder}} = \mathrm{mod}(\text{\texttt{work}}, N_p); \\
        &\text{\texttt{start}} = \text{\texttt{rank}} * \text{\texttt{chunk}} + \mathrm{min}(\texttt{remainder}, \text{\texttt{rank}}) + 1;\\
        &\text{\texttt{end}} = (\text{\texttt{rank}}+1) * \text{\texttt{chunk}} + \mathrm{min}(\text{\texttt{remainder}}, \text{\texttt{rank}}+1);\\
        &\text{\texttt{myChunk}} = \texttt{end} - \text{\texttt{start}}.
    \end{aligned}
\end{equation*}

At later iterations, when each process has less fine propagations to deal with, one can adjust using array indexing instead of dealing with multiple allocation/deallocation. Our strategy is presented in Algorithm \ref{AlgoMPI} and we shall now discuss the communication aspect of our approach.

\begin{algorithm}[ht]
    \caption{MPI implementation: multiscale kinetic parareal Algorithm \ref{para:algoparaBis}.}
    \label{AlgoMPI}
    \begin{algorithmic}[1]
    \Require $U^{0,0}$
    \For{$n = 1,\dots, N_g$} \Comment{First coarse guess}
        \State $U^{n,0} = \mathcal{G}\left(U^{n-1,0}\right)$ 
    \EndFor
    \vspace*{.2cm}
    \While{$k \leq K$ \textbf{or} error $\geq$ tol}\Comment{Parareal iterations}
    \vspace*{.2cm}
        \State \texttt{work} = $N_g-k$ \Comment{Work distribution}
        \vspace*{.05cm}
        \State \texttt{chunk} = $\left\lfloor\text{\texttt{work}} / N_p\right\rfloor$
        \vspace*{.05cm}
        \State \texttt{remainder} = $\mathrm{mod}(\texttt{work}, N_p)$
        \vspace*{.05cm}
        \State \texttt{start} = \texttt{rank} * \texttt{chunk} + $\mathrm{min}$(\texttt{remainder}, \texttt{rank}) + 1
        \vspace*{.05cm}
        \State \texttt{end} = (\texttt{rank}+1) * \texttt{chunk} + $\mathrm{min}$
        (\texttt{remainder}, \texttt{rank}+1)
        \vspace*{.05cm}
        \State \texttt{myChunk} = \texttt{end} - \texttt{start}
        \vspace*{.2cm}
        \For{$n = start,\dots, end$} \Comment{Compute the jumps in parallel}
            \State $\Delta^n = \mathcal{P}\mathcal{F}\mathcal{L}(U^{n-1,k-1}) - \mathcal{G}\left(U^{n-1,k-1}\right)$
        \EndFor
        \vspace*{.2cm}
        \State Gather chunks of jumps $\Delta^n$ to master process
        \vspace*{.2cm}
        \If{\texttt{rank} \textbf{is} master}
        \For{$n = k,\dots, N_g$} \Comment{Sequential correction}
            \State $U^{n,k+1} = \mathcal{G}\left(U^{n-1,k}\right) + \Delta^{n}$
        \EndFor
        \EndIf
        \vspace*{.2cm}
        \State Broadcast updated solution $U^{n,k+1}$;
        \State Compute successive error on the moments and 
        \State $k = k+1$
        \vspace*{.2cm}
    \EndWhile
    \end{algorithmic}
\end{algorithm} 

\medskip

\textbf{Communications.} The main challenge to implement this strategy will be to efficiently send time-chunks of data that are of size $5N_x\times\textsf{myChunk}$ which are potentially large. It is however important to note that since the workload (in time) diminishes with parareal iterations, an empirical threshold for the number of processors to use should be considered to keep the ammount of communications in check. Our approach is to gather all the local-in-time jumps on the master processor (step 14) and update the solution. Then, the full updated solution of size $5N_xN_g$ is broadcast to each processor (step 20) to perform its work for the next parareal iteration. Note that this strategy may not be the most efficient depending on the size of the problem and computing architecture. For example, we refer to \cite{ThakurRabenseifnerGropp2005} where message passing cost and optimization is discussed.

Finally, another computational difficulty will be the actual storage of such data that needs to be MPI friendly in order to avoid the use of many buffer arrays and have an optimized browsing of the memory.
    
\bibliographystyle{acm}
\bibliography{BiblioParareal}

\end{document}